\documentclass{article}
\usepackage{mitpress}
\usepackage{amsmath}
\usepackage{graphicx}
\usepackage{amsfonts}
\usepackage{amssymb}
\newtheorem{theorem}{Theorem}

\newtheorem{proposition}[theorem]{Proposition}
\newtheorem{remark}[theorem]{Remark}

\newdimen\dummy
\dummy=\oddsidemargin
\begin{document}

\title{Polyharmonic Daubechies type wavelets in Image Processing and Astronomy, I}
\author{Ognyan Kounchev, Damyan Kalaglarsky}
\maketitle

\textbf{Abstract: }\emph{We introduce a new family of multivariate wavelets
which are obtained by ''polyharmonic subdivision''. They generalize directly
the original compactly supported Daubechies wavelets.}


\section{Introduction}

We consider new \textbf{multivariate polyharmonic Daubechies type wavelets
}which are called ''polyharmonic subdivision wavelets''. They have been
recently introduced in the paper \cite{DynKounchevLevinRender}. They are
obtained by means of a procedure called ''polyharmonic subdivision'' which is
a generalization of the classical one-dimensional subdivision scheme of
Deslauriers-Dubuc \cite{deslaurierDubuc1989} which is the original source for
the first compactly supported wavelets of Daubechies in $1988$, cf.
\cite{daubechies}. This new family of polyharmonic wavelets is the second
representative of the \emph{Polyharmonic Wavelet Analysis} following the
''polyspline wavelets'' which have been introduced in the monograph
\cite{okbook}.

An important feature of these newly-born wavelets is that they are a nice
generalization of the one-dimensional wavelets of Daubechies: they form an
orthonormal family, enjoy nice non-stationary ''refinement operator''
equations, and have compact filters. In addition to that they have elongated
supports. Let us remind that a major drawback of the one-dimensional spline
wavelets of Ch. Chui is that they do not have finite filters, and
respectively, the polyspline wavelets of \cite{okbook} do not have finite filters.

\section{Construction of fundamental function $\Phi_{m}$ for exponential
polynomials subdivision}

The whole construction of the Daubechies type wavelets passes via the
construction of the so-called \emph{fundamental function of subdivision}, cf.
\cite{blatter}. In the present case we will work with \textbf{non-stationary}
subdivision and we have a family of such functions $\Phi_{m}$ for all
$m\in\mathbb{Z}$ which satisfy the \textbf{refinement equations} (two-scale
relations) given by
\begin{equation}
\Phi_{m}\left(  t\right)  =\sum_{i\in\mathbb{Z}}a_{i}^{\left[  m\right]  }%
\Phi_{m+1}\left(  2t-i\right)  \qquad\text{for all }t\in\mathbb{R}.
\label{refinement}%
\end{equation}
We define the \textbf{non-stationary subdivision symbol }by putting\textbf{ }
\begin{equation}
a^{\left[  k\right]  }\left(  z\right)  :=\sum_{j\in\mathbb{Z}}a_{j}^{\left[
k\right]  }z^{j}. \label{ak}%
\end{equation}

We are interested in special subdivision processes arising through the
solutions of Ordinary Differential Equations. We assume that we are given a
number of frequencies $0\leq\lambda_{1}\leq\lambda_{2}\leq...\leq\lambda_{p}$
and put for the frequency vector (with repetitions)
\[
\Lambda=\left\{  \lambda_{1},\lambda_{2},...,\lambda_{N}\right\}  \cup\left\{
-\lambda_{1},-\lambda_{2},...,-\lambda_{N}\right\}  .
\]
We consider the space of $C^{\infty}$ solutions of the ODE
\begin{equation}
\prod_{j=1}^{p}\left(  \frac{d^{2}}{dt^{2}}-\lambda_{j}^{2}\right)  f\left(
t\right)  =0. \label{Solutions}%
\end{equation}
Let us recall a simple fact from ODEs: in the case of different $\lambda_{j}%
$'s the space of all $C^{\infty}$ solutions in (\ref{Solutions}) is spanned by
the set $\left\{  e^{\lambda_{j}t}:j=1,2,...,p\right\}  .$ In the case of $s$
coinciding indices $\lambda_{i}=\lambda_{i+1}=...=\lambda_{i+s-1}$ we have
that the solution set contains the functions $\left\{  t^{\ell}e^{\lambda
_{i}t}:\ell=0,1,...,s\right\}  .$

Let us proceed to the construction of the subdivision symbols. We put
\[
x_{j}=e^{-\lambda_{j}/2^{k+1}}.
\]
We define the following Laurent polynomial
\[
d\left(  z\right)  :=d^{\left[  k\right]  }\left(  z\right)  :=\prod_{j=1}%
^{N}\frac{\left(  z+x_{j}\right)  \left(  z^{-1}+x_{j}\right)  }{\left(
1+x_{j}\right)  ^{2}}%
\]
and
\begin{equation}
P\left(  x\right)  :=P^{\left[  k\right]  }\left(  x\right)  :=\prod_{j=1}%
^{N}\left(  1-\frac{4x_{j}}{\left(  1+x_{j}\right)  ^{2}}x\right)  .
\label{Ppolynomial}%
\end{equation}
They satisfy the equality
\begin{equation}
d\left(  e^{i\omega}\right)  =P\left(  \sin^{2}\frac{\omega}{2}\right)
\qquad\text{for all }\omega\in\mathbb{R}; \label{d=P}%
\end{equation}
cf. \cite{micchelliWAnonstationary}. We will often drop the dependence on the
upper index in $d,$ $a,$ $P$ and the other functions and symbols.

An important step for construction of the subdivision coefficients
$a_{j}^{\left[  m\right]  }$ is the application of the Bezout theorem:

\begin{proposition}
\label{PBezoutMicchelli} There exists a unique polynomial $Q$ with real
coefficients of degree $N-1$ such that
\[
P\left(  x\right)  Q\left(  x\right)  +P\left(  1-x\right)  Q\left(
1-x\right)  =1
\]
and
\[
Q\left(  x\right)  >0\qquad\text{for }x\in\left(  0,1\right)  .
\]
\end{proposition}

We define now the trigonometric polynomial $b\left(  z\right)  =b^{\left[
k\right]  }\left(  z\right)  $ by putting
\[
b\left(  e^{i\omega}\right)  =Q\left(  \sin^{2}\frac{\omega}{2}\right)  .
\]
We finally define the symmetric Laurent polynomial $a\left(  z\right)  $ by
putting
\begin{equation}
a\left(  z\right)  :=a^{\left[  k\right]  }\left(  z\right)  :=2d\left(
z\right)  b\left(  z\right)  \qquad\text{for }z\in\mathbb{C}\setminus\left\{
0\right\}  . \label{aSymbol}%
\end{equation}

The following proposition is important for the application of the Riesz lemma
to $a\left(  z\right)  $ and construction of the Wavelet Analysis, cf.
\cite{micchelliWAnonstationary}, \cite{DynKounchevLevinRender}.

\begin{proposition}
\label{PMicchelli2} The polynomial $a\left(  z\right)  $ defined in
(\ref{aSymbol}) satisfies
\[
a\left(  z\right)  =\sum_{j=-2N+1}^{2N-1}a_{j}z^{j}%
\]
with $a_{j}=a_{-j}=\overline{a_{j}}$ and
\[
a\left(  z\right)  \geq0\qquad\text{for all }\left|  z\right|  =1.
\]
\end{proposition}

The following fundamental result shows that the symbols $a\left(  z\right)  $
are the non-stationary subdivision symbols for symmetric set of frequencies
$\Lambda,$ cf. \cite{DynKounchevLevinRender}.

\begin{theorem}
\label{Tsubdivision} For every exponential polynomial, i.e. for every solution
to the equation
\begin{equation}
Lf\left(  t\right)  :=\prod_{j=1}^{N}\left(  \frac{d^{2}}{dt^{2}}-\lambda
_{j}^{2}\right)  f\left(  t\right)  =0\label{Lf}%
\end{equation}
we put%
\[
f_{j}^{k}=f\left(  \frac{j}{2^{k}}\right)  .
\]
Then $f$ is reproduced by means of interpolatory subdivision, i.e.
\begin{align}
f_{j^{\prime}}^{k+1} &  =\sum_{j=-\infty}^{\infty}a_{j^{\prime}-2j}^{\left[
k\right]  }f_{j}^{k}\qquad\text{for all }j^{\prime}\in\mathbb{Z}%
\label{subdivision}\\
f_{2j}^{k+1} &  =f_{j}^{k}\qquad\qquad\ \qquad\text{for all }j\in
\mathbb{Z},\nonumber
\end{align}
For every $m\in\mathbb{Z}$ the fundamental function of subdivision $\Phi
_{m}\left(  t\right)  $ is a continuous function obtained throught the
subdivision process (\ref{subdivision}), where one starts from $f_{j}%
^{0}=\delta_{j}$ for $j\in\mathbb{Z}$ (here $\delta_{j}$ is the Kronecker
symbol), i.e. we put $\Phi_{m}\left(  \frac{j}{2^{m}}\right)  =\delta_{j},$
and $\Phi_{m}$ satisfies the refinement equation (\ref{refinement}).
\end{theorem}

Having in hand the functions $\Phi_{m}$ and their refinement symbols
$a^{\left[  m\right]  }$ we may follow the usual scheme for construction of
father and mother wavelets which has been used by Daubechies, cf.
\cite{daubechies}, \cite{blatter}. The following fundamental result has been
proved in \cite{DynKounchevLevinRender}.

\begin{theorem}
\label{Tbasic} There exists a polynomial $g\left(  z\right)  =\sum
_{j\in\mathbb{Z}}g_{j}z^{j}$ such that it is the ''square root'' of $2a\left(
z\right)  ,$ i.e.
\begin{equation}
a\left(  e^{i\theta}\right)  =\frac{1}{2}\left|  g\left(  e^{i\theta}\right)
\right|  ^{2} \label{a-squareroot}%
\end{equation}
For every $m\in\mathbb{Z}$ there exists a compactly supported function
$\varphi_{m}\left(  t\right)  $ which satisfies the refinement equation
\begin{equation}
\varphi_{m}\left(  t\right)  =\sum_{j}g_{j}\varphi_{m+1}\left(  2t-j\right)  ,
\label{firefinement}%
\end{equation}
and the family $\left\{  \varphi_{m}\left(  t-j\right)  \right\}
_{j\in\mathbb{Z}}$ is orthonormal. (These are the non-stationary father
wavelets.) The functions
\begin{equation}
\psi_{m}\left(  t\right)  =\sum_{j\in\mathbb{Z}}\left(  -1\right)  ^{j}%
g_{1-j}\varphi_{m+1}\left(  2t-j\right)  \label{psirefinement}%
\end{equation}
are the mother wavelets; the family $\left\{  \psi_{m}\left(  t-j\right)
\right\}  _{j\in\mathbb{Z}}$ is orthonormal and the family $\left\{  \psi
_{m}\left(  t-j\right)  \right\}  _{m,j\in\mathbb{Z}}$ forms an orthonormal
basis of $L_{2}\left(  \mathbb{R}\right)  .$
\end{theorem}

\subsection{The polyharmonic case}

For the polyharmonic subdivision we will work with very special ODEs defined
by $L_{\xi}:=\left(  d^{2}/dt^{2}-\xi^{2}\right)  ^{N}$ which are the Fourier
transform of the polyharmonic operator $\Delta^{N}.$ For a fixed constant
$\xi\geq0$ we put
\begin{equation}
\Lambda:=\left(  -\xi,-\xi,...,-\xi,\xi,\xi,...,\xi\right)  \in\mathbb{R}^{2N}
\label{Lambda}%
\end{equation}
i.e. $\lambda_{j}=\xi,$ for $j=1,2,...,N.$ Now for fixed $\xi\geq0$ and
$k\in\mathbb{Z}$ we define the polynomial
\begin{equation}
d\left(  z\right)  :=d^{\left[  k\right]  ,\xi}\left(  z\right)  :=d^{\left[
k\right]  }\left(  z\right)  :=\frac{\left(  z+x_{0}\right)  ^{N}\left(
z^{-1}+x_{0}\right)  ^{N}}{\left(  1+x_{0}\right)  ^{2N}}\qquad\text{for }%
z\in\mathbb{C}; \label{dPolynomial}%
\end{equation}
here we put $x_{0}:=e^{-\xi/2^{k+1}}.$ For the sake of simplicity we will very
often drop the dependence on $k$ and $\xi$. By (\ref{d=P}) we have $d\left(
e^{i\omega}\right)  =P\left(  \sin^{2}\frac{\omega}{2}\right)  $ where
\begin{equation}
P\left(  x\right)  =\left(  1-\frac{4x_{0}}{\left(  1+x_{0}\right)  ^{2}%
}x\right)  ^{N}=\left(  1-\eta x\right)  ^{N}, \label{PpolynomialPolyharmonic}%
\end{equation}
and we have put%
\[
\eta=\eta^{\left[  k\right]  ,\xi}:=\frac{4x_{0}}{\left(  1+x_{0}\right)
^{2}}=\frac{2}{1+\cosh\left(  \xi/2^{k+1}\right)  }.
\]
Then following Proposition \ref{PBezoutMicchelli} we have to find the
polynomial solution $Q$ to the equation
\[
P\left(  x\right)  Q\left(  x\right)  +Q\left(  1-x\right)  P\left(
1-x\right)  =1
\]
where $Q$ has degree $\leq N-1.$

\begin{remark}
Let us recall that the polynomial $Q$ in the classical case, cf. e.g.
\cite{blatter}, p. $195,$ satisfies condition
\[
\left(  1-y\right)  ^{N}Q\left(  y\right)  +y^{N}Q\left(  1-y\right)  =1.
\]
The lowest degree solution polynomial $Q$ will be called Daubechies'
polynomial and we put
\begin{equation}
R_{N}\left(  x\right)  :=\sum_{j=0}^{N-1}\binom{N+j-1}{j}y^{j}.
\label{DaubechiesPolynomial}%
\end{equation}
(Note that in \cite{daubechies} and \cite{blatter} the notation used is
$P_{N}$ !)
\end{remark}

It is \textbf{amazing} that it is possible to solve the problem in Proposition
\ref{PBezoutMicchelli} explicitly.

\begin{proposition}
\label{PQpolynomial} Let $\Lambda=\left(  -\xi,-\xi,...,-\xi,\xi,\xi
,...,\xi\right)  \in\mathbb{R}^{2N}.$ Then for the corresponding polynomial
$P\left(  x\right)  =\left(  1-\eta x\right)  ^{N},$ the polynomial $Q$ of
degree $N-1$ defined by
\begin{equation}
Q\left(  x\right)  =Q_{N}^{k,\xi}\left(  x\right)  =\left(  2-\eta\right)
^{-N}\sum_{j=0}^{N-1}\binom{N+j-1}{j}\frac{\left(  1-\eta\left(  1-x\right)
\right)  ^{j}}{\left(  2-\eta\right)  ^{j}} \label{Qpolynomial}%
\end{equation}
solves the equation
\begin{equation}
P\left(  x\right)  Q\left(  x\right)  +P\left(  1-x\right)  Q\left(
1-x\right)  =1. \label{Qequation}%
\end{equation}
Hence,
\begin{equation}
Q\left(  x\right)  =\left(  2-\eta\right)  ^{-N}R_{N}\left(  \frac
{1-\eta\left(  1-x\right)  }{2-\eta}\right)  . \label{Qpolynomial2}%
\end{equation}
\end{proposition}

Hence, we find the trigonometric polynomial $b^{\left[  k\right]  }\left(
z\right)  $ by putting
\begin{equation}
b^{\left[  k\right]  }\left(  z\right)  :=b^{\left[  k\right]  ,\xi}\left(
e^{i\omega}\right)  :=Q^{\left[  k\right]  ,\xi}\left(  \sin^{2}\frac{\omega
}{2}\right)  \label{bPolynomial}%
\end{equation}
where we recall the notations
\[
x=\sin^{2}\frac{\omega}{2}=\frac{1-\cos\omega}{2}=\frac{1}{2}-\frac{z+z^{-1}%
}{4},
\]
Finally, we obtain the subdivision symbol $a^{\left[  k\right]  }\left(
z\right)  $ by putting
\begin{equation}
a^{\left[  k\right]  }\left(  z\right)  :=a^{\left[  k\right]  ,\xi}\left(
z\right)  :=2d^{\left[  k\right]  ,\xi}\left(  z\right)  b^{\left[  k\right]
,\xi}\left(  z\right)  . \label{aPolynomial}%
\end{equation}
Now by Theorem \ref{Tbasic} we find the ''square root'' of the symbol
$a^{\left[  k\right]  }\left(  z\right)  .$ This means that we have to take
separately the ''square root'' of the Laurent polynomials $d^{\left[
k\right]  }\left(  z\right)  $ and $b^{\left[  k\right]  }\left(  z\right)  .$
The ''square root'' of $d^{\left[  k\right]  }\left(  z\right)  $ is obvious;
taking the ''square root'' of $b^{\left[  k\right]  }\left(  z\right)  $ needs
taking the ''square root'' of the polynomial $Q.$

\section{Algorithm for finding the square root of the polynomials $Q$}

For the algorithmic aspects of taking the ''square root'' of the polynomial
$Q$ it will be important to describe the polynomial $Q$ through the zeros of
the Daubechies' polynomial $R_{N}$ in (\ref{DaubechiesPolynomial}).

\begin{proposition}
\label{PDaubechiesZeros} Let the zeros of the Daubechies' polynomial
(\ref{DaubechiesPolynomial}) be $c_{j}^{D},$ i.e.
\[
R_{N}\left(  y\right)  =\sum_{j=0}^{N-1}\binom{N+j-1}{j}y^{j}=\frac{\left(
2N-2\right)  !}{\left(  \left(  N-1\right)  !\right)  ^{2}}\prod_{j=1}%
^{N-1}\left(  y-c_{j}^{D}\right)  .
\]
Then the polynomial $Q$ as determined by (\ref{Qpolynomial}) is given by
\[
Q\left(  x\right)  =\left(  2-\eta\right)  ^{-2N+1}\eta^{N-1}\frac{\left(
2N-2\right)  !}{\left(  \left(  N-1\right)  !\right)  ^{2}}\prod_{j=1}%
^{N-1}\left(  x-C_{j}\right)  ,
\]
where
\[
C_{j}:=\frac{c_{j}^{D}\left(  2-\eta\right)  +\eta-1}{\eta}.
\]
\end{proposition}

By formula (\ref{dPolynomial}) we have the representation
\[
d^{\left[  k\right]  }\left(  z\right)  =\left|  \frac{\left(  z+x_{0}\right)
^{N}}{\left(  1+x_{0}\right)  ^{N}}\right|  ^{2}\qquad\text{for }z=e^{i\omega
},
\]
hence, we take the trigonometric polynomial
\begin{equation}
M_{1}\left(  z\right)  :=\frac{\left(  z+x_{0}\right)  ^{N}}{\left(
1+x_{0}\right)  ^{N}} \label{M1}%
\end{equation}
as its ''square root'', i.e. $d^{\left[  k\right]  }\left(  z\right)  =\left|
M_{1}\left(  z\right)  \right|  ^{2}$ for $\left|  z\right|  =1.$ Further, we
have to take care of the ''square root'' of the polynomial $b^{\left[
k\right]  }\left(  z\right)  .$ Thus we have to find the polynomial $M_{2}$ of
degree $\leq N-1$ such that
\begin{equation}
\left|  M_{2}\left(  e^{i\omega}\right)  \right|  ^{2}=\frac{1}{2}Q\left(
\sin^{2}\frac{\omega}{2}\right)  , \label{M2square=Q}%
\end{equation}
which may be obtained by using the roots of the Daubechies polynomials.

\begin{remark}
Let the polynomial $Q$ have the zeros $C_{j}$ as in Proposition
\ref{PDaubechiesZeros}, and let us put
\[
c_{j}=1-2C_{j}.
\]
We see that $Q\left(  \sin^{2}\frac{\omega}{2}\right)  =\widetilde{Q}\left(
\cos\omega\right)  $ for some polynomial $\widetilde{Q}$ and $c_{j}$ are the
zeros of $\widetilde{Q}.$ Hence, we may apply the algorithm for the Riesz
representation of $\widetilde{Q}$, see e.g. \cite{blatter}, p. $197-198$.
\end{remark}

Thus we obtain finally for every integer $m\geq0$ and $\xi\in\mathbb{Z}^{n}$
the representation
\begin{equation}
a^{\left[  m\right]  ,\left|  \xi\right|  }\left(  z\right)  =\frac{1}%
{2}\left|  M_{1}\left(  z\right)  M_{2}\left(  z\right)  \right|  ^{2},
\label{a-factorization}%
\end{equation}
and the family of functions
\begin{equation}
M\left(  z\right)  :=M^{\left[  m\right]  }\left(  z\right)  :=M^{\left[
m\right]  ,\xi}\left(  z\right)  :=M_{1}\left(  z\right)  M_{2}\left(
z\right)  \label{Mm-refinement}%
\end{equation}
represents the refinement masks for the family of scaling functions (father
wavelets) $\left\{  \varphi_{m}\left(  t\right)  \right\}  _{m\geq0}$ for
which the functions $\Phi_{m}$ are autocorrelation functions.

\begin{remark}
Note that the above factorization has been found in the special case $\xi=0$
by Daubechies in \cite{daubechies}, p. $266$; the coefficients of the ''square
root'' polynomial for $N=2..10$ are in table $6.1$ in \cite{daubechies}. A
detailed discussion of more efficient methods for choosing the proper
polynomial $M_{2}\left(  z\right)  $ is available in Strang-Nguyen
\cite{strang}, p. $157,$ in chapter $5.4$ on Spectral factorization. The
factorization of the Daubechies' polynomial $R_{N}\left(  y\right)  $ is
discussed in Burrus \cite{burrus}, on p. $78$ and the Matlab program is
\textbf{[hn,hin]=daub(N) }in Appendix $C.$ They work with the zeros of the
polynomial $R_{N}$ and provide a number of manipulations for finding a more
stable factorization.
\end{remark}

\textbf{Acknowledgement.} The first named author was sponsored partially by
the Alexander von Humboldt Foundation, and both authors were sponsored by
Project DO--2-275/2008 ''Astroinformatics'' with Bulgarian NSF.

\textbf{ABOUT THE AUTHORS}

Ognyan Kounchev, Prof., Dr., Institute of Mathematics and Informatics,
Bulgarian Academy of Science, tel. $+359-2-9793851;$ kounchev@gmx.de

Damyan Kalaglarsky, Institute of Astronomy, Bulgarian Academy of Science, tel.
$+359-2-9793851;$ damyan@skyarchive.org.
\end{document}